\documentstyle[12pt,amssymb]{amsart}

\begin{document}

{\makeatletter
\gdef\eqalign#1{\null\,\vcenter{\openup\jot\m@th
  \ialign{\strut\hfil$\displaystyle{##}$&$\displaystyle{{}##}$\hfil
      \crcr#1\crcr}}\,}
}


\def\ale{\mathrel{\mathop{<}\limits_{\sim}}}
\def\age{\mathrel{\mathop{>}\limits_{\sim}}}
\def\gs{\mathrel{\mathop{>}\limits_{\sim}}}
\def\qed{{\bf q.e.d.}\par\medbreak}
\def\qedn{\thinspace\null\nobreak\hfill\hbox{\vbox{\kern-.2pt\hrule
height.2pt
depth.2pt\kern-.2pt\kern-.2pt \hbox to2.5mm{\kern-.2pt\vrule width.4pt
\kern-.2pt\raise2.5mm\vbox to.2pt{}\lower0pt\vtop to.2pt{}\hfil\kern-.2pt
\vrule width.4pt\kern-.2pt}\kern-.2pt\kern-.2pt\hrule height.2pt depth.2pt
\kern-.2pt}}\par\medbreak}
\def\pf{ \noindent{\sl Proof\/}.\enspace}
\let\emptyset\varnothing

\let\de=\partial
\def\eps{\varepsilon}
\def\phe{\varphi}
\def\Hol{\mathop{\mathrm{Hol}}\nolimits}
\def\Aut{\mathop{\mathrm{Aut}}\nolimits}
\def\Re{\mathop{\mathrm{Re}}\nolimits}
\def\Im{\mathop{\mathrm{Im}}\nolimits}
\def\Fix{\mathop{\mathrm{Fix}}\nolimits}
\def\id{\mathop{\mathrm{id}}\nolimits}
\def\cancel#1#2{\ooalign{$\hfil#1/\hfil$\crcr$#1#2$}}

\def\R{{\Bbb R}}
\def\P{{\Bbb P}}
\def\N{{\Bbb N}}
\def\Z{{\Bbb Z}}
\def\Q{{\Bbb Q}}
\def\C{{\Bbb C}}
\def\T{\hbox{\bf T}}

\newbox\bibliobox
\def\setref #1{\setbox\bibliobox=\hbox{[#1]\enspace}
        \parindent=\wd\bibliobox}
\def\biblap#1{\noindent\hang\rlap{[#1]\enspace}\indent\ignorespaces}
\def\art#1 #2: #3! #4! #5 #6 #7-#8 \par{\biblap{#1}#2: {\sl #3\/}.
        #4 {\bf #5}~(#6)\if.#7\else, \hbox{#7--#8}\fi.\par\smallskip}
\def\book#1 #2: #3! #4 \par{\biblap{#1}#2: {\rm #3.} #4.\par\smallskip}

\def\atanh{\mathop{\rm atanh}\nolimits}

\overfullrule=0pt
\let\smb=\smbar \let\j=\jmath \let\i=\imath
\lineskiplimit=1pt

\def\ca{Cara\-th\'{e}o\-dory~}
\def\ko{Ko\-ba\-ya\-shi~}
\def\psc{pseu\-do\-con\-vex~}
\def\psh{pluri\-sub\-har\-monic~}
\def\ngbd{neigh\-bor\-hood~}
\def\bdry{boundary~}
\def\O{\Omega}
\def\CC{{\cal C}}
\def\CH{{\cal H}}
\def\CN{{\cal N}}
\def\CO{{\cal O}}
\def\CS{{\cal S}}
\def\CF{{\cal F}}
\def\CZ{{\cal Z}}
\def\ls{\mathrel{\mathop{<}\limits_{\sim}}}
\def\cl{{\rm cl}}
\def\dist{{\rm dist}}
\def\Lip{{\rm Lip}}
\def\aut{{\rm Aut}}


\title{Finite Type Conditions on Reinhardt Domains}

\author[Siqi Fu, 
Alexander V. Isaev and
Steven G. Krantz]
{Siqi Fu\\ 
Alexander V. Isaev\\ 
Steven G. Krantz}

\address{\hskip-\parindent Siqi Fu,
Department of Mathematics,
University of California, 
Riverside, CA
92521-0135, USA}
\email{sfu@@math.ucr.edu}
\address{\hskip-\parindent A. V. Isaev,
Centre for Mathematics and Its Applications,
The Australian National University,
Canberra, ACT 0200,
Australia} 
\email{Alexander.Isaev@@anu.edu.au}

\address{\hskip-\parindent S. G. Krantz,
Department of Mathematics,
Washington University, St.~Louis, MO 63130,
USA }

\email{sk@@math.wustl.edu}

\thanks{Research supported in part by a grant
from the National Science Foundation.
Research at MSRI supported in part by NSF grant \#DMS 9022140.}

\maketitle

\begin{abstract} 
In this paper we 
prove that, if $p$ is a boundary point of a smoothly bounded 
pseudoconvex Reinhardt domain in $\C^n$, then the variety type at $p$ 
is identical to the regular type.
\end{abstract}

\vskip 30pt

In this paper we study the finite type conditions on 
pseudoconvex Reinhardt domain.  We prove that,
if $p$ is a boundary point of a smoothly bounded pseudoconvex Reinhardt
domain in $\C^n$, then the variety type at $p$ is identical 
to the regular type.  
In a forthcoming paper, we will study the biholomorphically invariant
objects
({\it e.g.}, the Bergman kernel and metric, the Kobayashi and
Carath\'{e}odory 
metrics) on a pseudoconvex Reinhardt domain of finite type.

\bigskip
We first recall some definitions.
A domain $\Omega\subset \C^n$ is Reinhardt if 
$$(e^{i\theta_1}z_1,\ldots,
e^{i\theta_n}z_n)\in \O$$ whenever $(z_1,\ldots, z_n)\in \O$  and
$0\le \theta_j\le 2\pi, 1\le j\le n$.   
Denote 
$$\CZ_j =\{ (z_1,\ldots, z_n)\in \C^n ;
\quad z_j=0\},$$ 
for $j = 1,\dots,n$.  Let $\CZ = \bigcup_{j=1}^{n} \CZ_j$.  Define
$L\colon \C^n\setminus \CZ \to \R^n$ by
$$
 L(z_1, \ldots, z_n) =(\log |z_1|, \ldots, \log
|z_n|) 
$$ 
and
$$
L^* (z_1, \ldots, z_n)=(\log z_1, \ldots, \log
z_n). 
$$
where, in the second case, the logarithm takes the principle branch, and
$L^*$ is defined locally near every point of $\C^n\setminus \CZ$. 



Let $\O=\{ z\in \C^n; \quad \rho (z) < 0\}$ be a bounded domain with
smooth boundary.  A boundary point $p\in b\O$ is of finite {\it variety
type m} if 
$$
m=\sup \left \{ {v(\rho\circ f)\over v(f)} \right \} <
\infty\leqno(1)  
$$
where the supremum is taken over all analytics disc $f\colon \Delta\to
\C^n$ such that $f(0)=p$.  We  use $v(f)$ to denote the order of
vanishing of $f$ at $0$ (see [K] for details on this notion).  
When the supremum in (1) is
taken over all  regular analytic discs $f$ ({\it i.e.},
$f'(0)\ne 0$),  then $m$ is called the {\it regular
type} of $p$; when the supremum in (1) is taken over all
complex lines through $p$, then $m$ is called the {\it
line type} of $p$.

It has been proved by McNeal [M] that, for a boundary point of a smoothly
bounded convex domain, the variety type is identical to the line type.
This result was generalized by Boas and Straube [B-S] to a star-like
boundary
point (see Theorem 4 below for the exact statement). 
Meanwhile, J. Yu [Y] showed that the Catlin multi-type is
identical to the D'Angelo q-type for a boundary point of a
convex domain.  In this paper, we shall prove that, for a
boundary point of a pseudoconvex Reinhardt domain, the
regular type is identical to the variety type.

\medskip
\proclaim{Lemma 1}.  (Reinhardt) If $\O$ is a pseudoconvex
Reinhardt domain, then $L(\O\setminus \CZ)$ is convex.  Furthermore,
if  $\O\cap \CZ_j\ne\emptyset$, then  $(z_1,\ldots, z_{j-1}, \lambda
z_j, z_{j+1}, z_n)\in\O$ whenever $(z_1, z_2, \ldots, z_n)\in\O$ and
$|\lambda|\le 1$. \par
\smallskip 
\proclaim{Lemma 2}. Let $\O$ be a smoothly bounded
pseudoconvex  Reinhardt domain and let $p\in b\O$.  Then
there is a neighborhood $U$ of $p$ such that $b\O\cap U$
has a plurisubharmonic defining  function $\rho (z)$.
\par   
\pf  If $p\in \O\setminus \CZ$,  then $L(z)$ is well-defined
in a neighborhood of $p$.  Since $L(\O\setminus \CZ)$ is convex,
there is a local convex defining function $\tilde\rho (u_1, \ldots,
u_n)$ of $L(\O\setminus \CZ)$ near $L(p)$.  Let 
$\rho =\tilde\rho (L(z))$.  Then $\rho$ is a plurisubharmonic
defining function of $b\O$ near $p$.

If $p\in \O\cap \CZ$,  then, without lose of generality,  we may assume
that
$$p=(p_1, \ldots, p_k, 0,\ldots, 0),$$  
with $p_j\ne 0$, $0\le j\le k$.
Let $\sigma (z)$  be the Euclidean distance from $z$ to $b\O$.
Then $\sigma (z)$ is smooth in a neighborhood of $b\O$ and satisfy
$\sigma (z_1, \ldots, z_n) =\sigma (e^{i\theta_1}z_1, \ldots, 
e^{i\theta_n}z_n)$. It follows that ${\de \sigma \over \de z_j}(z_1,
\ldots,
z_n) = 0$ whenever $z_j =0$. Since ${\de \sigma\over \de z_j}(p) =0$,
$k+1\le j \le n$, 
one of ${\de \sigma\over \de z_j} (p)$ , $1\le j\le k$, must 
not be zero.
 
For simplicity of notation, we may assume that 
${\de \sigma\over \de z_1} (p)\ne 0$. Therefore $b\O$ is locally
defined by  an equation of the form $|z_1| + f(|z_2|, \ldots, |z_n|)=0$
near
$p$.  This in turn implies that $L(b\O\setminus\CZ)$ is defined by
$$
u_1 + g(u_2, \ldots, u_n)=0
$$
where
$$
g(u_2, \ldots, u_n)=-\log(-f(e^{u_2}, e^{u_3}, \ldots,
e^{u_n}))  
$$
for $(u_2, \ldots u_n)\in I_2\times \cdots I_k\times I^{n-k}$.
Here $I_j$ is a neighborhood of $\log |p_j|$, $2\le j\le k$, and
$I=(-\infty, -C)$ ($C >0$ is a large constant).  

The convexity of $L(b\O\setminus \CZ)$ implies that $g(u_2, \ldots,
u_n)$ is a convex function.  Therefore
$$
\eqalign{\rho(z_1, \ldots, z_n)&= \log |z_1|+
                                g(\log |z_2|, \ldots, \log |z_n|)\cr
                               &=\log |z_1| - \log
(-f(|z_2|, \ldots, |z_n|))\cr
}                                      
$$
is a smooth plurisubharmonic defining function for $b\O$
near $p$. \qedn
\smallskip
is

formulate

\noindent For $z=(z_1, \ldots, z_n)$, let $\hat z =(z_2,\ldots, z_n)$.

\proclaim{Theorem 3}([B-S]). Suppose that a smooth
hypersurface is defined by $\Re z_1+h(\hat z, \Im z_1)=0$
near the origin.  If $h(\hat z, 0)\ge 0$, then the variety
type of the hypersurface
at the origin is obtained by analytic discs of form $\phi (\zeta)=
(0, \phi_2(\zeta), \ldots, \phi_{n}(\zeta))$.  Furthermore, if
the hypersurface is star-like at the origin, {\it i.e.},  
there exists a $\delta >0$ such that
$t\to h(a_2t, \ldots, a_{n}t, 0)$ is an increasing function on 
$[0, \delta ]$ for all unit vectors $(a_2, \ldots, a_{n})$,  then
the variety type at the origin is obtained by a complex line.
\par
Now we can prove the main theorem of this paper:

\proclaim{Theorem 4}. If $\O$ is a smoothly bounded
pseudoconvex Reinhardt domain and if $p\in b\O$, then the
variety type at $p$ is  identical to the regular type.
\par

\pf 
It follows from the smoothness and the Reinhardt
property of $b\O$ that $p$ cannot be the origin (this can be
easily seen  from the proof of Lemma 2 in this paper or
from Lemma 1.3 in [FIK]). We divide the proof into three cases.

If $p\in b\O\setminus \CZ$, then by Lemma 1, $L^*
(b\O\setminus\CZ)$ is convex near $p$. Therefore, it
follows from the results of McNeal [Mc] and Boas/Straube
[B-S] that the variety type of $L^* (b\O\setminus Z)$ at
$L^*(p)$ is obtained by a complex line of the form $L^*(p)
+a\zeta$. Thus, the variety type of $b\O$ at $p$ is
obtained by $$ \phi (\zeta) =(p_1e^{a_1 \zeta}, \ldots,
p_ne^{a_n\zeta}). $$
Therefore the variety type is identical to the regular type at $p$.
\smallskip\smallskip

If $p\in b\O\cap \CZ$, then the mapping $L^*$ is not well-defined near
$p$. Thus the above method fails. We first consider the case that 
only one coordinate of $p$ is non-zero.  We shall show that, in this 
case, the domain $\O$ is
star-like at $p$. After a simple dilation and
rotation, we may assume that $p=(1, 0,\ldots, 0)$ and $\Re z_1$  is the
outward normal direction at $p$.  It follows from the proof of
Lemma 2 that $b\O$ is defined near $p$ by 
$$
\log |z_1| + h(\hat z) =0
$$
where $h$ is a smooth function defined near $\hat p$ such that
\smallskip 

\begin{enumerate}
\item[(a)]$ h(\hat p) =0$ and $dh(\hat p)=0;$ 
\item[(b)] $h(e^{i\theta_2}z_2, \ldots, e^{i\theta_n}z_n)=h(z_2, \ldots,
            z_n);$
\item[(c)] The function $h$ is plurisubharmonic.  Furthermore, 
$$h(z_2,\ldots, z_n) = g(\log |z_2|, \ldots, \log |z_n|)$$
for some convex function $g$ and
$z_i\ne 0$, $2\le i\le n$.
\end{enumerate}

\smallskip\smallskip
For $a=(a_2,\ldots, a_n)$ with $\| a\| =1$ and $a_i\ne 0$, $2\le i \le
n$, 
define
$$
\eqalign{ k_a(t)&=h(ta_2, \ldots, ta_n)\cr
               &=g(\log t + \log |a_2|, \ldots, \log t +\log |a_n|)
.\cr} 
$$
It follows from the convexity of $g$ that $\tilde{g}(s)=g(s +\log |a_2|,
\ldots, s+\log |a_n|)$ is a convex function of $s$. By (a), we have
$$
\lim_{s\to -\infty}\tilde{g}'(s)=0. 
$$
Therefore, $\tilde{g}'(s) \ge 0$ for $s<<-1$. It then follows that
$k_a(t)$ is an increasing function of $t$. Thus, by the
(proof of the) theorem of Boas/Straube, 
the variety type at $p$ is obtained by a complex line of the form
$\phi(\zeta) = (1, b_2 \zeta, \ldots, b_n\zeta).$  

\smallskip
Now we treat the general situation for the case when $p\in b\O\cap
\CZ$. Without lost of generality, we may assume that the first $k$
coordinates of $p$ are 1's and the rest of the coordinates are 0's. 
Let $\sigma (z)$ be the  Euclidean distance from $z$ to $b\O$. 
It follows from the proof of Lemma 2 that one of ${\de \sigma \over
\de z_j}(p)$, $1\le j \le k$ is not zero. For simplicity of notation, we
assume that  $\partial \sigma (p)/\partial z_1 \ne 0$.  Let $z'=\Phi (z)$
be defined by
$$
\eqalign{z'_1& =z^{\alpha_1}_1\cdots z^{\alpha_k}_k;\cr
          z'_j&=z_j, \quad 2\le j\le n\cr}
$$
where $\alpha_j =\partial \sigma(p)/\partial z_j, 1\le j\le k$.

It is obvious that $z'=\Phi(z)$ is biholomorphic in a neigborhood
of $p$ and $\Re z'_1$ is the outward normal direction to $b\O$ 
at $p$.  This change of coordinates preserves the Reinhardt
property of the
domain near $p$.  Changing coordinates if necessary, we may 
assume at the outset that $\Re z_1$ is the outward 
normal direction to
$b\O$ at $p$.  It follows from the proof of Lemma 2 that
$b\O$ is defined near $p$ by
$$
\log |z_1| + h(z_2, \ldots, z_n) = 0 ,
$$
where $h$ is a function satisfying properties (a)-(c).  By a similar
argument, we can prove that, for fixed $z_j$, $2\le j\le k$, near 1 and
$a=(a_{k+1},\ldots, a_n)$ with $\|a\|=1$, 
$$
t \mapsto h(z_2,\ldots, z_k, ta)
$$
is an increasing function of $t$ for small $t>0$.  Thus we have
$$
h(z_2, \ldots, z_n)\ge h(z_2, \ldots, z_k, 0, \ldots,
0).\leqno(2) 
$$

Let $\rho(w_2,\ldots, w_n)= h(e^{w_2}, \ldots, e^{w_k}, w_{k+1},
\ldots, w_n)$. Then $\rho(\hat w)$ is a plurisubharmonic function. 
Furthermore, for fixed $(w_{k+1}, \ldots, w_n)$, $\rho(\hat w)$ is
a convex function.
Let $k(t)=\rho(tw_2, \ldots, tw_{k}, 0, \ldots, 0)$.  Then $k(t)$ is
a convex function with $k(0)=k'(0)=0$.  Therefore $k(1)=\rho(w_2, \ldots,
w_k, 0,\ldots, 0)\ge 0$. Hence $h(\hat z)\ge 0$.  
\smallskip

\proclaim{Claim}. The variety type of $b\O$ at $p$ is
obtained by analytic discs of the form $\phi(\zeta)=(1, \hat \phi(\zeta))$.
\par
The proof of the claim comes directly from the  proof of the theorem of
Boas and Straube [B-S]. We provide details here for the reader's
convenience. 

Let $r(z)=\log |z_1| + h(z_2, \ldots, z_n)$. Let 
$ \phi(\zeta)=(\phi_1(\zeta), \hat \phi(\zeta))\colon \Delta\to \C^n$
be an analytic disc such that $\phi(0)=p$. We may assume that 
$v(r\circ \phi)> v(\phi)$. If $\phi_1(\zeta)$ is not identically
zero, then there exists a sequence $\{ \zeta_j \}\subset \Delta$,
$\zeta_j\to 0$
such that $\{ \phi_1(\zeta_j) \}$ is a sequence of real numbers 
and $\phi_1(\zeta_j)\to 
1^+$.  Therefore 
$$
\log |\phi_1(\zeta_j)|=|\log \phi_1(\zeta_j)| \le  
r\circ \phi (\zeta_j).
$$ 
Let $\psi(\zeta)=(1, \hat \phi(\zeta))$.
Since the order of vanishing of a holomorphic 
function is determined by its order of vanishing along a sequence, we
obtain
that
$$
v(\phi_1(\zeta)-1)=v(\log |\phi_1(\zeta)|)\ge v(r\circ \phi).
$$ 
This in turn
implies that $v(\psi)=v(\phi)$.  It then follows from the equality
$$
|r\circ\psi(\zeta) -r\circ\phi(\zeta) |=|\log |\phi_1(\zeta)| |
$$
that $v(r\circ\psi)\ge v(r\circ \phi)$. This finishes the proof of
the claim.   
\smallskip\smallskip

We now return to the proof of the main theorem. 
Let $\tilde \phi(\zeta) =(\phi_2(\zeta),\ldots, \phi_k(\zeta))$ and
let $m=v(h\circ \hat\phi)$. It follows from (2) that
$$
v(h(\tilde \phi(\zeta),0, \ldots, 0))\ge m.\leqno(3)
$$
Let $\beta_j =v(\phi_j)$, $k+1\le j\le n$ and $\psi_j(\zeta)=\phi_j(\zeta)
/\zeta^{\beta_j}$.  For $\delta>0$, we have
$$
\eqalign{
\delta^m&\gs {1\over 2\pi\delta}\int_{|\zeta|=\delta} 
h\circ \hat\phi(\zeta) |d\zeta|\cr
        &={1\over 2\pi\delta}\int_{|\zeta|=\delta} 
       h(\tilde\phi(\zeta),\delta^{\beta_{k+1}}\psi_{k+1}(\zeta), 
           \ldots \delta^{\beta_n}\psi_n(\zeta))|d\zeta|\cr
       &\ge h(1, \ldots, 1, 
              \delta^{\beta_{k+1}}\psi_{k+1}(0), \ldots 
                  \delta^{\beta_n}\psi_n(0)).\cr   
}
$$

Therefore,
$$
v(h(1, \ldots, 1, \zeta^{\beta_{k+1}}\psi_{k+1}(0), \ldots 
\zeta^{\beta_n}\psi_n(0)))\ge m.\leqno(4)
$$ 
  
It follows from (3) and (4) that the variety type at $p$
can be realized by analytic discs of the form 
$$
\phi(\zeta)=(1, \phi_2(\zeta),\ldots, \phi_k(\zeta), 0,\ldots, 0) 
\quad {\rm or}\quad
\phi(\zeta)=(1,\ldots, 1, \phi_{k+1}(\zeta),\ldots, \phi_n(\zeta)) .
$$
However, it follows from the first and second cases of this proof
that the variety type can actually be realized by analytic discs of 
form
$$
\phi(\zeta)=(1, e^{a_2\zeta}, \ldots, e^{a_k\zeta}, 0, \ldots, 0),\leqno(5)
$$
or
$$
\phi(\zeta)=(1, \ldots, 1, a_{k+1}\zeta, \ldots, a_n\zeta).\leqno(6)
$$

Indeed, to obtain the special form (5), consider the section of $\O$ by
$\CZ_{k+1,\dots,n}=\bigcap_{j=k+1}^n\CZ_j$. Then, by Lemma 1.3 of [FIK], this
section is a nonempty smoothly bounded set which is a finite collection of
Reinhardt domains in $\CZ_{k+1,\dots,n}$. The point $p$ then becomes a
point with all nonzero components, and therefore the first case of the
proof applies.

To obtain the special form (6), one needs to consider the section of $\O$
by a complex affine subspace of the form
$\{z_{j_1}=\dots=z_{j_{k-1}}=1\}$, where $1\le j_i\le k$ for all $i$. Then
we
can use the argument from the second case of the proof.  
\qedn

\noindent{\bf Remark 5.}\enspace  
{(1)} One certainly
 cannot expect that the line
type at $p$ is equal to the variety type. For example, let $\O$
be a pseudoconvex Reinhardt domain in $\C^2$ such that $b\O$ is 
locally defined by
$$
\rho(z_1, z_2)=\log |z_1| +\log |z_2| +(\log |z_1|-\log |z_2|)^4
$$
for $z=(z_1, z_2)$ near $p=(1, 1)$.  Then the line type at $p$ is 2
while the regular type at $p$ is 4.

\noindent{(2)} The result of Yu [Y] cannot be generalized to
Reinhardt domains. 
Let $\Omega$ be a pseudoconvex Reinhardt domain defined by
$$
\rho(z_1, z_2, z_3)=|z_1|^2 + |z_2|^6 + |z_3|^6 + |z_2 z_3|^2 - 1 .
$$
Then the 3-tuple of the D'Angelo q-types at $p=(1, 0, 0)$ is
$(\Delta_3, \Delta_2, \Delta_1)=(1, 4, 6)$ while the Catlin 
multi-type is $(1, 4, 4)$. We thank Professor Yu for pointing
out this fact.

\end{document}